\documentclass[12pt]{amsart}
\usepackage{amsmath}
\usepackage{amsthm}
\usepackage{amscd}

\newtheorem{lemma}{Lemma}%               with section number.  Same
\newtheorem{corollary}{Corollary}%       goes for lemmas, etc.)

\theoremstyle{definition}

\theoremstyle{remark}

\begin{document}

\title[Archimedean analogue of Tate's conjecture]{On an archimedean
analogue of Tate's conjecture}

\author{Dipendra Prasad}
\email{dprasad@mri.ernet.in}

\author{C. S. Rajan}
\address{Harish-Chandra  Research Institute, 
Chhatnag Road, Jhusi, Allahabad 211019, 
INDIA}

\address{Tata Institute of Fundamental 
Research, Homi Bhabha Road, Bombay - 400 005, INDIA.}
\email{dprasad@mri.ernet.in, rajan@math.tifr.res.in}


{\sf

\begin{abstract}
We consider an Archimedean analogue of Tate's conjecture,  
and  verify the conjecture in the examples of isospectral
Riemann surfaces constructed by Vign\'eras and Sunada.
We prove a simple lemma in
group theory which lies at the heart of T. Sunada's theorem about
isospectral manifolds.

\end{abstract}

\maketitle

The aim of this short note is to formulate
a conjecture about isospectral Riemann
surfaces, which can be considered as an Archimedean analogue of Tate's
conjecture, and to  verify the conjecture in the examples of isospectral
Riemann surfaces constructed by Vign\'eras and Sunada. 
 We first  prove a simple lemma in
group theory which lies at the heart of T. Sunada's theorem about
isospectral manifolds (i.e., Riemannian manifolds for which the
eigenvalues of the Laplacian, counted with multiplicity, are the same).
The lemma has many other applications to similar problems. 

We begin with the group theoretic lemma. Let $G$ be a finite
group.  For a  $G$-module $M$, denote by  $M^G$
the submodule of invariants of $G$. 

\begin{lemma}
Suppose that $G$ is a finite group with subgroups $H_1$ and $H_2$ such
that each conjugacy class in $G$ intersects $H_1$ and $H_2$ in equal
number of elements. 
Assume that $V$ is a representation space of $G$ 
over a field $k$  of characteristic zero. Then there exists an isomorphism
$i: V^{H_1} \rightarrow V^{H_2}$, 
commuting with the action of any endomorphism 
$\Delta$ of $V$ which commutes with the action of $G$ on $V$, i.e. the 
following diagram commutes:

$$
\begin{CD}
  {\rm V}^{H_1} @>i>> {\rm V}^{H_2} \\
  @V{\Delta}VV    @VV{\Delta}V \\
  {\rm V}^{H_1} @>i>> {\rm V}^{H_2}
\end{CD}
$$

Further if $k= {\Bbb C}$ and 
the action of $G$ on $V$ is unitary, the isomorphism can be
chosen to be unitary. 
\end{lemma}

\begin{proof} The hypothesis on the subgroups $H_1$ and $H_2$ 
of $G$ implies that the character of the left regular representation of $G$ 
on the space of functions on $G/{H_1}$ with values in $k$, 
$k[G/{H_1}]$, is the same as the character of the representation
of $G$ on $k[G/{H_2}]$. Hence these representations are 
isomorphic (over $k$!). The isomorphism is well-known to exist by character 
theory over an algebraically closed field containing $k$, and hence the two
representations are isomorphic over $k$ too by a general result in group
representations: two representations of a group $G$ over 
a field $k$ which 
become isomorphic over a field extension of $k$, are isomorphic over $k$; cf. 
page 110 of the article of Atiyah and Wall in the book edited by J.W.S.
Cassels and A. Fr\"ohlich on `Algebraic Number Theory'. The article of Atiyah 
and Wall deals with $k={\Bbb Q}$ which is all that is needed by us, however the same
proof works for any infinite field (and can be proved for finite fields too 
by a simple application of Lang's theorem).

By Frobenius reciprocity, for any   
representation $V$ of $G$, there is an isomorphism of
 $ V^H$ with $[V \otimes {k}[G/H]]^G$ 
in which we send a vector 
$v \in V^H$
 to
$ \sum gv \otimes e_g$, 
the sum running over the distinct elements of $G/H$, 
 denoted by $e_g$. 

The isomorphism of  $ V^H$ with $[V \otimes {k}[G/H]]^G$ for $H=H_1,H_2$,
taken together with a
 $G$-equivariant isomorphism 
$S: k[G/{H_1}] \rightarrow k[G/{H_2}],$ gives an isomorphism of $V^{H_1}$ with
$V^{H_2}$:
$$V^{H_1}  \rightarrow [V \otimes {k}[G/H_1]]^G \rightarrow 
[V \otimes {k}[G/H_2]]^G \rightarrow V^{H_2}.$$

Any endomorphism of $V$ 
 which commutes with the action of $G$
on $V$ preserves the invariant subspaces $V^{H_1}$ and $V^{H_2}$. 
It is clear that the isomorphism between $V^{H_1}$ and $V^{H_2}$ constructed 
here commutes with such endomorphisms. 

It remains to check the unitarity. Let $V$ be a unitary $G$-module. Since 
${\Bbb C}[G/H_1]$ and  ${\Bbb C}[G/H_2]$ are isomorphic as
$G$-modules, we can choose a unitary  isomorphism with the natural
unitary structures on  ${\Bbb C}[G/H_1]$ and ${\Bbb C}[G/H_2]$. This
is a general fact: if there are two unitary structures on a complex 
representation $V$ of a group $G$, then there exists a $G$-invariant
intertwining operator between the two unitary structures; 
we omit the proof here.
For each $g\in G$, the elements $ gv \otimes e_g$ are
 mutually orthogonal, each one of  norm equal to norm of $v$, where
$e_g$ denotes the characteristic function of the coset space $gH$ in
$G/H$. 
Hence it follows that the isomorphism of
 $ V^H$ with $[V \otimes {\Bbb C}[G/H]]^G$ in which we send a vector 
$v \in V^H$ to
$ \frac{1}{\sqrt{|G/H|}}\sum_{g\in G/H} gv \otimes e_g$ is unitary, 
if $V \otimes {\Bbb C}[G/H]$ 
is given the usual  unitary structure  for a tensor product. Thus 
if $V$ is a unitary $G$-module, we obtain  a unitary  isomorphism from
$V^{H_1}$ to $V^{H_2}$.

\end{proof}

\begin{corollary}[Sunada's theorem] 
Suppose that $G$ is a finite group with subgroups $H_1$ and $H_2$ such
that each conjugacy class in $G$ intersects $H_1$ and $H_2$ in equal
number of elements. Suppose that $X$ is a Riemannian manifold on which 
$G$ acts by isometries and freely (i.e., $gx=x$ implies 
$g=e$). Then the quotient of $X$ by $H_1$ and $H_2$ with the induced
metrics are isospectral manifolds. 
\end{corollary} 

\begin{proof} Let $V$ be the space of ${\mathcal C}^{\infty}$ 
functions on $X$.
This is a representation of $G$ such that $V^{H_i}$ is the space of 
${\mathcal C}^{\infty}$ functions on the quotient of $X$ by $H_i$. The
Laplacian on $X$ commutes with the $G$-action, and when restricted
to functions in $V^{H_i}$, it corresponds to the Laplacian on
${\mathcal C}^{\infty}$ functions on the quotient of $X$ by $H_i$. 
The proof of the corollary hence follows from the Lemma. 
\end{proof}
 
\begin{corollary} 
Let $H$ be a Lie group, $\Gamma_0, \Gamma_1, \Gamma_2, \Gamma$
discrete subgroups of $G$ with 
$$\Gamma_0 \subset \Gamma_i \subset \Gamma, ~~~~~~{\rm~~~~~~~for~~~~~~ }i=1,2$$
such that $\Gamma_0$ is a normal subgroup of finite index in $\Gamma$ 
with the property that the subgroups $\Gamma_i/ \Gamma_0$ of $\Gamma/\Gamma_0$
for $i=1,2$ intersect each conjugacy class in $\Gamma/\Gamma_0$ in equal 
number of elements. Then $L^2(\Gamma_1\backslash H)$ and
$L^2(\Gamma_2\backslash H)$ are isomorphic as $H$-modules.
\end{corollary}

\begin{proof}  The proof of the corollary follows from Lemma 1 applied 
to $V = L^2(\Gamma_0\backslash H)$,  with $G = 
\Gamma_0\backslash \Gamma$, $H_1=\Gamma_0\backslash \Gamma_1$ and 
$H_2 = \Gamma_0\backslash \Gamma_2$. We note here that a similar corollary
is also available in the work of Vign\'eras, corollary 5 of [V].
\end{proof}

\begin{corollary} Let $X$ be a projective algebraic variety over a field $k$
together with an action of a finite group $G$ on $X$ over $k$. Let 
$\overline{k}$ be a separable closure of $k$, and $\overline{X}$ denote 
$X$ base changed to
$\overline{k}$.   
Let $H_1$ and $H_2$ be two subgroups of $G$ which intersect each conjugacy
class in $G$ in equal number of elements. 
Then for any $i$,
$H^i_{\text{{\'e}t}}(\overline{X},{\Bbb Q}_{\ell} )^{H_1}$ and 
$H^i_{\text{\'{e}t}}(\overline{X},{\Bbb Q}_{\ell} )^{H_2}$ 
are isomorphic as Gal$(\overline{k}/k)$-modules. Equivalently,
for any $i$,
$H^i_{\text{\'{e}t}}(\overline{X}/{H_1},{\Bbb Q}_{\ell} )$ and 
$H^i_{\text{\'{e}t}}
(\overline{X}/{H_2},{\Bbb Q}_{\ell} 
)$ are isomorphic as Gal$(\overline{k}/k)$-modules. Hence for $k$ a finite
field, or a number field, the $L$-functions associated to 
$H^i_{\text{\'{e}t}}(\overline{X}/{H_1},{\Bbb Q}_{\ell} )$ and 
$H^i_{\text{\'{e}t}}
(\overline{X}/{H_2},{\Bbb Q}_{\ell} 
)$ are the same. 
\end{corollary}

\begin{corollary}\label{jacobians}
Let $X$ be a projective algebraic curve over a field $k$ with
an action of a finite group $G$ on $X$ over $k$. 
Let $H_1$ and $H_2$ be two subgroups of $G$ which intersect each conjugacy
class in $G$ in equal number of elements. Then the Jacobians of the
curves $X/{H_1}$ and $X/{H_2}$ are isogenous over $k$.
\end{corollary}
\begin{proof}
 It can be seen, cf. lemma below, that 
the Jacobian $J_i$ of $X/{H_i}$ is
isogenous to the connected component of the $H_i$ fixed points 
of the Jacobian $J$ of $X$. Hence once again,
the map $$J^{H_1} \rightarrow J^{H_2}$$
defined by $x \rightarrow \sum \Phi(g)g\cdot x$ provides the necessary
isogeny from $J_1$ to $J_2$. (Here $\Phi$ is an
integral valued function on $G$ corresponding to an isomorphism of
${\Bbb Q}[G/H_1]$ with ${\Bbb Q}[G/H_2]$, interpreted as a 
function on the double coset space $H_2\backslash G/H_1$ .)
\end{proof}

The following lemma is well-known. However, not finding an appropriate
reference, we have included a proof here.

\begin{lemma}
Let $X$ be a projective algebraic curve over a field $k$ with
an action of a finite group $H$ on $X$ over $k$. Then
the Jacobian of $X/{H}$ is
isogenous to the connected component of the $H$ fixed points 
of the Jacobian $J$ of $X$. 
\end{lemma}
\begin{proof} Let $n$ be the order of $H$. 
We can assume that $k$ is algebraically closed. 
We will prove that multiplication by $n^2$ takes the group of $H$-fixed 
points on $J$ to the natural image of the Jacobian of 
${X}/H$ into the
Jacobian of ${X}$, proving the lemma.
Let ${k}(X)$ denote the function field of ${X}$, 
and ${\rm Div}^0({X})$,
the divisors of degree 0 on ${X}$. 
We have the exact sequence of $H$-modules,
$$0\rightarrow {k}(X)^*/{k}^* \rightarrow 
{\rm Div}^0({X}) \rightarrow 
{\rm Pic}^0({X})\rightarrow 0.$$
Since  $H^1(H,A)$ is annihilated by $n$ for any $H$-module $A$, 
we have that any $H$-invariant element in 
${\rm Pic}^0({X})$,  when multiplied by $n$ comes from an
$H$-invariant divisor on ${X}$.  
But  any $H$-invariant divisor on ${X}$ multiplied by  $n$ comes
from ${X}/H$,  completing the proof of the lemma.
\end{proof}

\noindent{\bf Remark 1.}
Many examples of triples $(G;H_1,H_2)$ as above have been known for a 
very long time.  Some of these are given in Sunada's paper; see also 
Serre's book, {\it Linear Representations of finite groups}, section 
13.2, exercises 5 and 6.

\vspace{4mm}

\noindent{\bf Remark 2.}
 Corollary 3 specialises to give examples of (non-conjugate) 
number fields with the same zeta function. The extension of this
phenomenon to geometric context was Sunada's original motivation 
for his paper. 

\vspace{4mm}

Corollary 4 lends
support to the following conjecture, which can be considered as an
Archimedean analogue of a corollary of 
Tate's conjecture and proved by Faltings in [F], 
that if the Hasse-Weil zeta
functions of two curves defined over a number field $K$ are equal, then
the Weil restriction of scalars from $K$ to ${\Bbb Q}$ 
of the corresponding Jacobians are isogenous. 

\vspace{4mm}

\noindent{\bf Conjecture.} {\it Suppose that $X$ and $Y$ are complete algebraic curves 
defined over a number field $K$. Suppose that the compact Riemann
surfaces associated to $X$ and $Y$, via an embedding of 
$K$ into ${\Bbb C}$, are isospectral with respect to the K\"ahler metric of
constant curvature (1, 0, or $-1$ as the case may be). 
Then the Hasse-Weil zeta functions of $X$ and $Y$
are the same over a finite extension $L$ of $K$, and hence the Weil 
restriction of scalars from $L$ to ${\Bbb Q}$ 
of the corresponding Jacobians are isogenous.}

\vspace{4mm}

\noindent{\bf Remark 3.} Let $X$ and
$X'$ be two Shimura curves defined using quaternion division algebras 
over totally real number field. It is known by a result of A. W. Reid
(Duke J., volume 65 (1992)), that if the Shimura curves are 
isospectral then the underlying quaternion algebras
are isomorphic.

\vspace{4mm}

We briefly recall the construction of isospectral surfaces by Vign\'eras
\cite{V}, and verify the above conjecture.
 Let $K$ be a totally real number field, and let ${\Bbb H}$
be a  
quaternion division algebra over $K$ which is ramified at all but one of
the real places $v_0$.  Suppose that ${\mathcal O}$ and ${\mathcal O}'$ are two
maximal orders in ${\Bbb H}$, which are not 
conjugate by any ${\Bbb Q}$-automorphism of ${\Bbb H}$. 
 Let $\Gamma$ (resp. $\Gamma'$) be
the  group of elements of reduced norm $1$ in ${\mathcal O}$
(resp. ${\mathcal O}'$) modulo the group ${\pm 1}$. Projecting to the
$v_0$-component,   $\Gamma $ and $\Gamma'$
give rise to  co-compact lattices in $PSL(2,{\Bbb R})$. 
Let $X$ and $X'$ be the
corresponding Riemann surfaces.  Under further
technical conditions Vign\'eras shows that $X$ and 
$X'$  are isospectral but not
isometric.  
By theorem 2.5 of Shimura \cite{Sh}, the curves $X$ and
$X'$ are defined over a certain number   field $K$, and that $X$
is isomorphic to a Galois conjugate of $X'$. This implies in particular 
that the Hasse-Weil zeta functions of $X$ and $X'$ are the same over 
some number field.

\vspace{4mm}

\noindent{\it Note added in Proof:} We do not know if the
following much stronger form of the above conjecture is true:
If two Riemann surfaces  (not necessarily defined over a number field) 
are isospectral, then the
Jacobian of one is isogenous to a conjugate of the other 
by an automorphism 
$\sigma \in {\rm Aut}({\Bbb C}/{\Bbb Q})$, where
$\sigma$ preserves the spectrum of the corresponding 
Riemann surface. More specifically, for 
Riemann surfaces defined over ${\Bbb Q}$, we do not know if isospectral
implies isogeny of the Jacobians.

\vspace{4mm}

\noindent{\bf Acknowledgement:} The authors would like to thank 
the referee for a very careful reading, and for pointing out a serious
error in an earlier version of the paper.


\begin{thebibliography}{99}
{\sf

\bibitem[F]{F} G. Faltings, {\it Endlichkeitss\"{a}tze f\"{u}r
    abelsche Variet\"{a}ten \"{u}ber Zahlk\"{o}rpern}, Invent. Math.,
    73 (1983) 349-366.

\bibitem[Sh]{Sh} G. Shimura, {\it On canonical models of arithmetic
quotients of bounded symmetric domains}, Annals of Mathematics, 91, (1970),
    144-222. 

 
\bibitem[S]{S}
T. Sunada, {\it Riemannian coverings and isospectral manifolds}, Annals 
of Mathematics, 121, (1985), 169-186.

\bibitem[V]{V} M.-F.Vign\'eras, {\it Vari\'et\'es Riemanniennes 
isospectrales et non
isom\'etriques}, Annals of Mathematics, 112, (1980), 21-32.
}
\end{thebibliography}
\end{document}